\documentclass[12pt]{amsart}

\usepackage{amscd}
\usepackage{amssymb,amsmath,amsthm}
\theoremstyle{plain}
\newtheorem{theorem}{Theorem}[section]
\newtheorem{corollary}[theorem]{Corollary}

\newtheorem{proposition}[theorem]{Proposition}

\theoremstyle{definition}
\newtheorem{definition}[theorem]{Definition}
\newtheorem{remark}[theorem]{Remark}
\newtheorem{ex}[theorem]{Example}
\def\ses#1#2#3{ $$
      \begin{CD}
      0 @>>> #1 @>>>#2 @>>> #3 @>>> 0
      \end{CD}
       $$  }

\def\mu #1 { \mathcal{M}(#1)}

\def\tsr{\operatorname{tsr}}

\begin{document}
\title [Some examples of lifting problems from quotient algebras ]
       {Some examples of lifting problems from quotient algebras}

\begin{abstract}
 We consider three lifting questions: Given a $C\sp{*}$-algebra $I$,
 if there is a unital $C\sp{*}$-algebra $A$ contains $I$ as an ideal,
 is every unitary from $A/I$ lifted to a unitary in $A$? is
 every unitary from $A/I$ lifted to an extremal partial
 isometry? is every extremal partial isometry from $A/I$
 lifted to an extremal partial isometry? We show several
 constructions of $I$
 which serve as working examples or counter-examples for above questions.
\end{abstract}

\author { Hyun Ho \quad Lee }
\address {Department of Mathematics\\
         Purdue University\\
         West Lafayette, Indiana 47907 }
\email{ylee@math.purdue.edu}
\keywords{Extremally rich
$C\sp{*}$-algebras, Good index theory, Stable rank one }
\subjclass[2000]{Primary:46L05.}
\date{Sep 15, 2008}
\maketitle

\section{Introduction}
Let $\mathcal{E}(A)$, or just $\mathcal{E}$, denote the set of
 extreme points of the unit ball of a unital $C\sp{*}$-algebra $A$. Recall that
 elements in $\mathcal{E}$ are characterized as the partial
 isometries $u$ in $A$ satisfying
 $(1-uu^{*})A(1-u^{*}u)=0$ by R. V. Kadison \cite{ka}. We call them extremal partial
 isometries and call the projections $1-uu^{*}$, $1-u^{*}u$ defect projections.
 In \cite{brpe2}, Brown and Pedersen defined the notion of extremal
richness for $C\sp{*}$-algebra $A$ which means quasi-invertible
elements are dense in $A$ as an analogue of stable rank one for
possibly infinite $C\sp{*}$-algebras. (We say $T$ in $A$ is
quasi-invertible if $T$ has closed range and the kernel projections
of $T^{*}$ and $T$ are centrally orthogonal in $A$, or if $T$ is in
$A^{-1}\mathcal{E}A^{-1}$. For more equivalent definitions, see
Theorem 1.1 in \cite{brpe2}.) We denote by $A_{q}^{-1}$  the set of
 quasi-invertible elements. As a result, stable rank one
$C\sp{*}$-algebras are characterized within the class of extremally
rich $C\sp{*}$-algebras by the property that all extreme points of
the unit ball are unitaries or $A_{q}^{-1}=A^{-1}$ where $A^{-1}$ is the set of invertible elements of $A$.\\
Suppose we have an extension of $C\sp{*}$-algebras: \ses{I}{A}{B} It
is well known that if we have an extension of an extremally rich
$C\sp{*}$-algebra $I$ by an extremally rich $C\sp{*}$-algebra $B$,
we cannot deduce $A$ is extremally rich even in the finite case. The
obstacle, as in the analogous problem with stable rank one, can be
expressed as a lifting problem but with special properties. In fact,
Brown and Pedersen proved the following theorem \cite[Theorem
6.1]{brpe}.

\begin{theorem}\label{T:E}
If $A$ is an extension of $I$ by $B$ and both $I$ and $B$ are
extremally rich, then $A$ is extremally rich if and only if extremal
partial isometries in $B$ are liftable to extremal partial
isometries in $A$ and $PAQ$ is extremally rich whenever $P$ is
projection
 of the form $1-uu^{*}$ where $u \in \mathcal{E}(A)$ and $Q$ projection
of the from $1-v^{*}v$ where $v \in \mathcal{E}(\tilde{I})$.
\end{theorem}

 However, there are some special ideals $I$  for which the hypotheses in Theorem \ref{T:E}
 can be simplified:
 \begin{corollary}
 Let $I$ be a $C\sp{*}$-algebra of stable rank one. Then $A$ is extremally rich
 if and only if $B$ is extremally rich and extremal partial isometries in $B$
 lift.
 \end{corollary}
 This corollary implies the following: Whenever $I$ embedded as an
 ideal in a unital extremally rich $C\sp{*}$-algebra and $u$ is an extremal partial isometry in
 a extremally rich $C\sp{*}$-algebra $A/I$, there is an extremal
 partial isometry in $A$ which lifts $u$.

 Motivated by above example, but not necessarily related to extension, we can ask whether we can find examples
 of $I$ such that whenever there exists $A$ containing $I$ as a closed ideal,  certain
 lifting questions below from any quotient $C\sp{*}$-algebra $B=A/I$ to $A$ are affirmative:
\begin{enumerate}
 \item[(1)] Is every unitary in $B$ lifted to a unitary in $A$?
 \item[(2)] Is every unitary in $B$ lifted to an extremal partial isometry in $A$?
 \item[(3)] Is every extremal partial isometry in $B$ lifted to an extremal partial
 isometry in $A$?
\end{enumerate}

We cannot expect a positive answer in great generality due to the
nature of the problem; In fact, it is not very difficult to find
counter examples for these questions. Interesting direction,
however, is to find an affirmative answer. For example, it is
well-known that if $I$ has an approximate identity consisting of
projections, then every unitary in $B$ lifts to a partial isometry
in $A$.

\begin{remark}\label{R:ext}

\begin{itemize}
\item[(i)] When $I$ and $B$ are of stable rank one, an affirmative answer
to (1) is equivalent to the stable rank one property for $A$.
\item[(ii)] When $I$ and $B$ are of stable rank one, an affirmative answer
to (2) is equivalent to extremal richness for $A$.
\item[(iii)] When $I$ has stable rank one and $B$ is extremally rich, an
affirmative answer to (3) is equivalent to extremal richness for
$A$.
\end{itemize}
\end{remark}

\section{Examples}
Throughout this article $H$ will denote  an infinite dimensional
Hilbert space and $B(H)$ the set of bounded operators on $H$ and
$K(H)$(shortly $K$) the set of compact operators on $H$. Given a
$C\sp{*}$-algebra $I$, $M(I)$ will denote the multiplier algebra of
$I$ and $C(I)$ will denote the corona algebra of $I$, that is,
$M(I)/I$.
 In view of Remark \ref{R:ext}-(i), it is still interesting if we restrict $I$ to be the class of stable rank one
 $C\sp{*}$-algebras.  But, even in this case,  question (1) is not true for every stable
 rank one $C\sp{*}$-algebras. In fact, there is a well-known counter example: the Toeplitz extension \ses
{K}{\mathcal{T}}{C(S^{1})} where $\mathcal{T}$ is the
$C\sp{*}$-algebra generated by the unilateral shift $S$ on $H$ is
such an example by the remark \ref{R:ext}-(i) and the fact
$\tsr(\mathcal{T}) \ne1$.(See \cite[4.13]{rie}.) However, there is a
$C\sp{*}$-algebra of stable rank one which serve as an answer for
question (1).

\begin{definition}(\cite[p.2]{brpe3})
  We say  a (non-unital) $C\sp{*}$-algebra $I$ has good index
 theory
 if whenever $I$ is embedded as an ideal in a unital $C\sp{*}$-algebra $A$
 and $u$ is a unitary in $A/I$ such that $\partial_{1}([u])=0 $ where
 $\partial_{1}:K_{1}(A/I) \to K_{0}(I)$, there is a unitary in $A$
 which lifts $u$.
\end{definition}
It was noted by Brown and Pedersen that G. Nagy proved that any
stable rank one $C\sp{*}$-algebra has good index theory
\cite{brpe3}. Thus if we can show the existence of $I$ with
$\tsr(I)=1$ such that $\partial_1:K_1(B) \to K_0(I)$ is trivial,
question (1) holds for $I$.
\begin {ex}
Let $$I=\{(f_{n})\mid f_{n} \in C(\mathbb{T},M_{n}(\mathbb{C}))\,
\mbox{and} \, f_{n} \to f \,\mbox{in} \,C(\mathbb{T},\mathcal{K})
\}$$  Here $f_{n} \to f \,\mbox{in} \,C(\mathbb{T},\mathcal{K})$
means $ \sup_{t\in \mathbb{T}}\|f_{n}(t)-f(t) \|$ goes to $0$ as $n$
goes to infinity.
 Note that $I$ has stable rank one.  Let $C(\mathbb{T},B(H)_{*-S})$
 be the set of functions $m:\mathbb{T} \to B(H)$ which are
 continuous with respect to the $*$-strong operator topology on
 $B(H)$, and which satisfy $\|m\|_{\infty}:=\sup_{t\in \mathbb{T}}\|m(t)\| <
 \infty$.

 It is not hard to show that
 $$ M(I)=\{(F_{n}) \mid F_{n} \in C(\mathbb{T},M_{n}(\mathbb{C
})) \,\mbox{and} \, F_{n}\to_s F  \, \mbox{in} \,
C(\mathbb{T},B(H)_{*-S}) \}$$ Here $ F_{n} \to_s F \, $ means
$F_nf_n \to Ff $ and $F_n^*f_n \to F^*f$ in $I$ for each $(f_n) \in
I$. Next we want to show $K_{0}(\iota):K_{0}(I)\to K_{0}(M(I))$ is
injective where $\iota: I \to M(I)$; since the sequence
$\alpha=(z_{n})$ in $K_{0}(I)$ is eventually constant, we may
   assume $I$ has only finite dimensional irreducible representations. Let
   $\pi$ be such an irreducible representation and suppose
   $K_{0}(\iota)(\alpha)=0$. Then $K_{0}(\bar{\pi})(K_{0}(\iota)(\alpha))=0 $
   where $\bar{\pi}:M(I) \to B(H)$ is a unique map such that
   $\bar{\pi}\circ \iota =\pi$. Note that we have $\pi(I)=\bar{\pi}(I)$ since
   $\pi$ is finite dimensional representation. Consequently, we have
   $K_{0}(\pi)(\alpha)=0$ in $K_{0}(\pi(I))$. By the Remark 5.12 in
   \cite{brpe}, $\alpha=0$.\\
   From the diagram
   \begin{equation*}
   \begin{CD}
   0 @>>> I @>j>>A @>>> B @>>> 0 \\
    @. @VVV @VVV @VVV \\
  0 @>>> I @>>\iota>M(I) @>>> C(I) @>>> 0
   \end{CD}
   \end{equation*}
   we have
    \begin{equation*}
   \begin{CD}
   K_{0}(I) @>K_0(j)>> K_{0}(A) \\
     @|  @VVV \\
   K_{0}(I) @>>K_0(\iota)> K_0(M(I))
   \end{CD}
   \end{equation*}
   Thus it follows that $K_0(j)$ is also injective.
   Thus from the exactness of six term
   exact sequence of K-theory we know the map $\partial_1:K_{1}(B) \to K_0(I)$ is trivial.
   \end{ex}
 Next we show there are  counter-examples to question (2).
 \begin{ex}\label{E:example2-1}
    Let $$ I=\{ (a_{n}) \in \prod M_{2n}(\mathbb{C})\mid
    a_{n} \to \left(
                        \begin{smallmatrix}
             D & 0 \\
             0 & E
             \end{smallmatrix}
               \right) \,\text{in} \quad M_{2}(K)\}$$
    Then
             $$M(I)=\{ (T_{n}) \in \prod M_{2n}(\mathbb{C})\mid
             T_{n} \to \left(
                        \begin{smallmatrix}
             M & 0 \\
             0 & N
             \end{smallmatrix}
               \right) \,\text{in} \quad M_{2}(B(H)_{*-s}) \}$$
               Now let's fix the basis $\mathfrak{B}$ of Hilbert space $H$ as follows;
                $$
               \{\cdots,w_{n},\cdots,w_{2},w_{1},v_{1},v_{2},\cdots
               v_{n},\cdots \} $$
               Let $e_{x \, y}$ be a rank one projection such that
               $e_{x \, y}(z):=<y,z>x$, and take
               $T_{n}=\sum_{i=1}^{n-1} e_{w_{i} \, w_{i+1}} +
               \sum_{j=1}^{n-1} e_{v_{j+1} \, v_{j}}+
                e_{w_{n} \, v_{n}}$. In fact, $T_{n}$ is of the following form with respect to the basis
               $\mathfrak{B}$:
                 \begin{equation*}
                  \begin{pmatrix}
    0     & 0      &  \dots &0      &0      & \dots &  0     &1\\
    1     & 0      &  \dots &0      &0      & \dots &  0     &0\\
    \vdots& \ddots & \vdots &\vdots &\vdots &\vdots &\vdots  &\vdots\\
    0     & 0      &  1     &0      &0      & \dots &0       &0\\
    0     & \dots  &  \dots &0      &0      & \dots &0       &0\\
    0     & 0      &  \dots &0      &1      & \dots &0       &0\\
    \vdots& \vdots & \ddots &0      &0      & \ddots&0       &0\\
    0     & \dots  & \dots  &0      &0      & \dots &1       &0
              \end{pmatrix}
              \end{equation*}
    Then $1-T_{n}^{*}T_{n}=e_{w_{1}\,w_{1}}$,
   $1-T_{n}T_{n}^{*}=e_{v_{1}\,v_{1}}$. Therefore $\pi((T_{n}))$ is an
   unitary element of $C(I)$ where $\pi $ is the natural quotient map.
   But it cannot be liftable to an extremal partial isometry: Assume there is an
   element $(a_{n}) \in I$ such that $T_{n}+a_{n}$ is extremal partial isometry
   in $M(I)$. Since $ M_{2n}(\mathbb{C})$ has stable rank one,
   $T_{n}+a_{n}$ must be a unitary, and it happens only when
   $a_{n}=e_{v_{1}\,w_{1}}$. Then $ T_{n}+a_{n}$ cannot converge to the
   operator of the form $\left(
                        \begin{smallmatrix}
             * & 0 \\
             0 & *
             \end{smallmatrix}
               \right)$. It is a contradiction.
 \end{ex}
More sophisticated example can be found as follows.
\begin{ex}\label{E:example2-2}
Let $D$ be the stabilization of cone or suspension. Assume there
exist a projection $\mathbf{p}$ in the corona algebra of $D$ such
that which does not lift but its $K_0$-class does lift.(For the
construction of such a projection, see the example 5.13 in
\cite{lee}.) If we let $a$ be the self adjoint element which lifts
$\mathbf{p}$ in $M(D)$, we take $I$ to be the $C\sp{*}$-algebra
generated by $a$ so that the quotient $I/D$ is isomorphic to
$\mathbb{C}$. Then the Busby invariant is determined by sending $1$
to $\mathbf{p}$, and we have the following commutative diagram:
\begin{equation*}
   \begin{CD}
   0 @>>> D @>j>> I @>>> \mathbb{C} @>>> 0 \\
    @. @| @VVV @VV\mathbf{p}V \\
  0 @>>> D @>>>M(D) @>>> C(D) @>>> 0
   \end{CD}
   \end{equation*}
  By the long exact sequence, we have
\begin{equation*}
   \begin{CD}
    0 @>>>K_0(I) @>>> K_0(\mathbb{C}) @>\partial_0>> K_1(D)  \\
     @. @VVV   @VV\mathbf{p}V @| \\
   0 @>>>K_0(M(D)) @>>> K_0(C(D)) @>\partial_0>> K_1(D)
   \end{CD}
   \end{equation*}
   Since $\partial_0([\mathbf{p}]_0)=0$, $\partial_0:K_0(\mathbb{C}) \to
   K_1(D)$ becomes trivial. Thus $K_0(I) \cong K_0(\mathbb{C})$. In
   particular, $K_0(I)$ is non-trivial. Consequently, we found a
   (stably) projectionless stable rank one $C\sp{*}$-algebra such that
   its $K_0$-group is non-trivial. Now stabilize this algebra and
   call it $I$ again. We consider an extension of $I$ by
   $C(\mathbb{T})$ corresponding to a unitary $u$ in the coronal algebra
   with non-trivial class $K_1$-class. $u$ can't be lifted to a
   unitary (If so, $[u]_1=0$ which is a contradiction), and it can't
   be lifted to a partial isometry either because there are no
   non-zero projections available to defect projections of the
   partial isometry.
\end{ex}

On the other hand, we have a certain class of
$C\sp{*}$-algebras(e.g. elementary $C\sp{*}$-algebras) such that
question (2) is true.
\begin{proposition}\label{P:question2}
Let $I$ be the $C\sp{*}$-algebra such that $M(I)$ is extremally
rich. Suppose a unital $C\sp{*}$-algebra $A$ contains $I$ as an
ideal. Then any unitary $u$ in $A/I$ can be liftable to an extremal
partial isometry.
\end{proposition}
\begin{proof}
From the diagram
\begin{equation*}
   \begin{CD}
   0 @>>> I @>>>A @>>> A/I @>>> 0 \\
    @. @VVV @VVV @VV\tau V \\
  0 @>>> I @>>>M(I) @>>> C(I) @>>> 0
   \end{CD}
\end{equation*}
we see $\tau(u)$ is also a unitary in $C(I)$. Since $M(I)$ is an
extremally rich $C\sp{*}$-algebra, we can find an extremal partial
isometry $w$ which is an inverse image of $\tau(u)$ under the map
$\pi$ so that we have $\pi(w)=\tau(u)$. Hence $(w,u)$ is the
extremal partial isometry in $A$ which lifts $u$ since we can view
$A$ as the pullback construction of $A/I$ and $M(I)$.
\end{proof}
\begin{remark}
We do not know the result of Proposition \ref{P:question2} is true
under the weaker hypothesis $C(I)$ is  extremally rich. In addition,
$I$ in Example \ref{E:example2-1} and Example \ref{E:example2-2} are
(stable rank one) $C\sp{*}$-algebras such that $M(I)$ is not
extremally rich.
\end{remark}

\begin{corollary}\label{C:question3}
Let $I$ be the $C\sp{*}$-algebra such that $M(I)$ is extremally rich
and any Busby invariant into $C(I)$ is extreme point preserving map.
Suppose a unital $C\sp{*}$-algebra $A$ contains $I$ as an ideal.
Then any extremal partial isometry $u$ in $A/I$ can be lifted to an
extremal partial isometry.
\end{corollary}
\begin{proof}
If a map is extreme point preserving, it sends an extremal partial
isometry to an extremal partial isometry. Thus the result follows
from above diagram.
\end{proof}

We have seen an example of $I$ such that any extremal partial
isometry in $A/I$ can't be lifted to a partial isometry in $A$ for
some $C\sp{*}$-algebra $A$ which contains $I$ as an ideal in Example
\ref{E:example2-2}. If this does not happen, we can show the lifting
question (3) might be true for a certain class of idealizers of $I$.
\begin{proposition}
Let $I$ be the $C\sp{*}$-algebra such that any extremal partial
isometry in $A/I$ is lifted to a partial isometry whenever a
$C\sp{*}$-algebra $A$ contains $I$ as an ideal. If $A$ is an
extremally rich $C\sp{*}$-algebra, then any extremal partial
isometry in $A/I$ can be lifted to an extremal partial isometry in
$A$.
\end{proposition}
\begin{proof}
Let $q$ be the quotient map from $A$ to $A/I$. If $u$ is an extremal
partial isometry in $B=A/I$, it can be lifted to a partial isometry
$v$ in $A$. Since $A$ is an extremally rich $C\sp{*}$-algebra, $v$
has an extremal extension $w$ in $\mathcal{E}(A)$ such that
$v=wv^*v=vv^*w$ \cite[Proposition 2.6]{brpe2}. Now let $I^{+}$ and
$I^{-}$ be the defect ideals of $u$ and consider the map
$\pi^{+}\oplus \pi^{-}: B \to B/I^{+}\oplus B/I^{-}$ which is
injective. Since
$\pi^{+}\oplus\pi^{-}(q(v))=\pi^{+}\oplus\pi^{-}(q(w))$,
$q(v)=q(w)=u$. Hence an extremal partial isometry $w$ is an inverse
image of the extremal partial isometry $u$.
\end{proof}
However, the following examples
   show us that affirmative examples for question (3) might be hard to get without putting a condition on
   $I$.
   \begin{ex}
   Let \begin{equation*}
    A= \left\{ \begin{pmatrix}
             A & B \\
             C & D
             \end{pmatrix}
                 \in M_{2}(B(H)) \biggm|  B,C \in
               K \right\}
         \end{equation*}
Then $I=M_{2}(K)$ is an ideal of $A$. And
    $B=A/M_{2}(K) $ is isomorphic to $ Q(H)\oplus Q(H)$.
    Now consider $ \begin{pmatrix}
             S & 0 \\
             0 & T^{*}
             \end{pmatrix}$ where $S$ and $T$ are isometries such that
               $1-SS^{*}$ and $1-TT^{*}$ are not compacts. Then $\pi \left(
            \begin{pmatrix}
             S & 0 \\
             0 & T^{*}
             \end{pmatrix}
               \right)$ becomes an extremal partial isometry in $Q(H)\oplus
               Q(H)$ where $\pi:I
               \to B$. We can check, however, that there is no extremal partial isometry which
               lifts  $\pi \left(
                      \begin{pmatrix}
             S & 0 \\
             0 & T^{*}
             \end{pmatrix}
               \right)$: Suppose there is $R \in
               M_{2}(K) $ such that $V=
                      \begin{pmatrix}
             S & 0 \\
             0 & T^{*}
             \end{pmatrix}
                + R$ is an extremal partial isometry. Since
               $A$ contains $M_{2}(K)$ as an ideal, either one
               of defect projections of $V$ must be zero. But note that $1-V^*V=
                      \begin{pmatrix}
             0 & 0 \\
             0 & 1-TT^{*}
             \end{pmatrix}
               + \mbox{compact}$ is in $M_{2}(K)$.
               Therefore  $1-VV^*=0$ implies $1-TT^*$ is compact which
               is a contradiction. The other case is also similar.
\end{ex}

   The following clever example is due to Larry Brown although it
   was slightly generalized by the author.
   \begin{ex}
   Let $A$  be a unital $C\sp{*}$-algebra such that $RR(A)=0$, $\tsr(A)=1$ and
   $K_{1}(A)=0$. Then it is known that $\tsr(C[0,1]\otimes A)=1$ \cite{mop}. Now let's
   consider $I$ as $C[0,1]\otimes (A \otimes \mathcal{K})$ which has stable
   rank one too. Note that $M(A\otimes \mathcal{K})_s$ is equipped with the strict
   topology. Then $ M(I) = C([0,1], M(A\otimes K)_s ) $.\\
    If we denote $[\mathcal{F}_{A}]$ by the set of homotype classes
    of Freehold elements in $ M (A\otimes K) $, there is Mingo's index
   map $\partial :[\mathcal{F}_{A}] \to K_0(A)$  which is actually an
   isomorphism \cite{mi}. Suppose $K_0(A)$ is an ordered group and let $e$ be an order
   unit. Then let $f_0 \in M(I)$ such that $f_0(1)=1$ and $f_0(t)$ is an
   isometry with $\partial(f_0(t))=-2e$ for $t<1$. Let $u$ be a co-isometry in
   $M(A\otimes K)$ with the index $e$ and $f$ be $f_0 u$. Then we can
   see that $\pi(f)$ is an isometry in $\mathcal{C}(I)$ where $\pi:M(I)\to C(I)$.
   Assume  there is a $k \in I$ such that $f+ k$ is an extremal partial isometry in $M(I)$
   which is a continuous path of extremal partial isometries.
Since  $M(A \otimes K)$ is a prime $C\sp{*}$-algebra, using index
theory, $f(1) + k (1)$ must be co-isometry of index $e$. It follows
that there is a unitary $v$ such that $u + k(1)=uv$ \cite[Theorem
2.1]{brpe4}. Since $uv-u$ is in $A \otimes K$, $u$ is a Fredholm
element in $M(A\otimes K)$, and $v \in 1 + A\otimes K$.\\
 Finally, let $f_1= f v$ and $k_1$ be $f+k-f_1=f(1-v)+k \in I$. Note
 that $k_1(1)=0$ and $f+ k=f_1+ k_1$. Thus $f_1(t)+ k_1(t)$ must be
an isometry for $t<1$ since $f$ has negative index for $t<1$ and
$f(t)+k(t)$ is an extremal partial isometry (isometry or co-isometry
in this case) for each $t$. Since $f_1$ has a non-trivial kernel
which is the range of $v^*(1-u^*u)$, we know that $|| k(t) || \geq
1$ for $t < 1$ which is a contradiction.
\end{ex}
  \section{Acknowledgements}
  This is a part of the author's Ph.D. thesis. He would like to
  thank his advisor Larry Brown for many helpful suggestions and
  discussions.

   \end{document}